\documentclass[reqno,11pt]{amsart}

\usepackage[top=2.0cm,bottom=2.0cm,left=3cm,right=3cm]{geometry}
\usepackage{amsthm,amsmath,amssymb,dsfont}
\usepackage{mathrsfs,amsfonts,functan,extarrows,mathtools}
\usepackage[colorlinks]{hyperref}

\usepackage{marginnote}
\usepackage{xcolor}
\usepackage{stmaryrd}
\usepackage{esint}
\usepackage{graphicx}
\usepackage{bm}
\usepackage{soul}
\newtheorem{theorem}{Theorem}[section]
\newtheorem{proposition}{Proposition}[section]

\newtheorem{remark}{Remark}[section]

\numberwithin{equation}{section}
\allowdisplaybreaks

\def\R{\mathbb{R}}

\def\T{\mathbb{T}}

\def\v{\varepsilon}

\def\l{\langle}
\def\r{\rangle}

\def\S{\mathbb{S}}

\def\up{\textup}
\def\n{\nabla}
\def\p{\partial}

\def\a{\alpha}
\def\b{\beta}

\def\c{\cdot}

\newcounter{wronumber}

\begin{document}
\title[Formal derivations from Boltzmann equation to three stationary equations]
			{Formal derivations from Boltzmann equation to three stationary equations}

\author[Zhendong Fang]{Zhendong Fang}
\address[Zhendong Fang]
        {\newline School of Mathematics, South China University of Technology, Guangzhou, 510641, P. R. China}
\email{zdfang@scut.edu.cn}

\thanks{2023}

\maketitle

\begin{abstract}
In this paper, we concentrate on the connection between Boltzmann equation and stationary equations. To our knowledge, the stationary Navier-Stokes-Fourier system, the stationary Euler equations and the stationary Stokes equations are formally derived by moment estimate in the first time and extend the results of Bardos, Golse, and Levermore in J. Statist. Phys. 63(1-2), 323-344, 1991.\\

\noindent {MSC:} 35B25;35Q30; 35Q20 \\

	\noindent{Keywords:} Hydrodynamic limit, Moment estimate, The stationary Navier-Stokes-Fourier system, The stationary Euler equations, The stationary Stokes equations.

\end{abstract}





\section{Introduction}
\subsection{The Boltzmann equation}
The Boltzmann equation is the fundamental equation in kinetic theory that describes the motion of molecules in phase space, which is written in the following form:
\begin{equation}\label{BE equaion}
\p_t f+v\c\n_xf=Q(f,f)
\end{equation}
where $f(t,x,v)$ denotes the number density of the gas molecules at time $t\geq0$, with position $x\in\Omega$ ($\Omega$ is $\T^3$ or $\R^3$) and velocity $v\in\R^3$, the collision operator $Q(f,f)$ is given by
\begin{equation}
Q(f,g)=\int_{\R^3}\int_{\S^2}(f'g'_*-fg_*)b(|v-v_*|,\omega)d\omega dv_*,
\end{equation}
here $b(|v-v_*|,\omega)$ denote collision kernel, $f=f(t,x,v),\,g_*=g(t,x,v_*),\,f'=f(t,x,v'),\,g'_*=g(t,x,v'_*)$, $(v,v_*)$ are velocities of two particles before collision while $(v',v'_*)$ are velocities of two particles after collision, $\omega$ is one of the (external) bissectors of the angle between two vectors $v-v_*$ and $v'-v'_*$ which means $|\omega|=1$, $(v',v'_*)$ are given by
\begin{equation}\label{Out veocities}
\begin{aligned}
v'=v'(v,v_*,\omega)=&v-[(v-v_*)\c\omega]\omega,\\
v'_*=v'_*(v,v_*,\omega)=&v_*+[(v-v_*)\c\omega]\omega.
\end{aligned}
\end{equation}

The Boltzmann equation, as a kinetic theory equation, has a deep connection with the fluid equations. After nondimensionalization\cite{SRL09}, we can deduce the scaled Boltzmann equation by introducing two dimensionless quantities of the kinetic $Strouhal$ number St, the inverse $Knudsen$ number Kn
\begin{equation}\label{Scaled BE equaion-0}
\up{St}\p_t f+v\c\n_xf=\frac{1}{\up{Kn}}Q(f,f).
\end{equation}

This is a multi-scale equation about St and Kn. We choose $\up{St}=\v,\,\up{Kn}=\v^q\,(q>0)$, then the scaled Boltzmann equation is rewritten as
\begin{equation}\label{Scaled BE equaion-1}
\v\p_t f_\v+v\c\n_xf_\v=\frac{1}{\v^q}Q(f_\v,f_\v).
\end{equation}

Formally, suppose $f_\v\to f$ as $\v\to 0^+$, using \eqref{Scaled BE equaion-1}, then
\begin{equation}\label{Formally-0}
\v^{q+1}\p_t f_\v+\v^q v\c\n_xf_\v=Q(f_\v,f_\v),
\end{equation}
which imply that
\begin{equation}\label{Formally-1}
Q(f_\v,f_\v)=\v^{q+1}\p_t f_\v+\v^q v\c\n_xf_\v\to 0,
\end{equation}
meanwhile, we have
\begin{equation}\label{Formally-2}
Q(f_\v,f_\v)\to Q(f,f).
\end{equation}

Therefore, combining with \eqref{Formally-1} and \eqref{Formally-2},we have
\begin{equation}\label{Formally-3}
Q(f,f)=0.
\end{equation}

According to the Boltzmann's $H$ theorem \cite{GFSRL05}, the following conditions are equivalent:
\begin{itemize}
    \item $Q(f,f)=0$ a.e.,
    \item $\int_{\R^3}Q(f,f)ln f dv=0$,
    \item $f$ is a Maxwellian density, i.e.
    \begin{equation}
     f=\mathcal{M}_{(\rho,u,\theta)}:=\frac{\rho}{(2\pi\theta)^{\frac{3}{2}}}\up{e}^{-\frac{|v-u|^2}{2\theta}}
    \end{equation}
    for some $\rho,\theta>0$ and $u\in\R^3$.
\end{itemize}

By applying Galilean transformation, we seek a special form solution near $M:=\mathcal{M}_{(1,0,1)}$ for the scaled Boltzmann equation \eqref{Scaled BE equaion-1}
\begin{equation}\label{Special form solution}
f_\v=M+\v^r g_\v M
\end{equation}
for some $r>0$.

Replacing \eqref{Special form solution} with the scaled Boltzmann equation \eqref{Scaled BE equaion-1} to get
\begin{equation}\label{Scaled BE equation}
\v\p_t g_\v+v\c\n_x g_\v+\frac{1}{\v^q}L g_\v=\frac{1}{\v^{q-r}}\Gamma(g_\v,g_\v),
\end{equation}
where the linearized Boltzmann operator $L$ and the bilinear symmetric operator $\Gamma(g,g)$ are given by
\begin{equation}\label{Linear opeators}
\begin{aligned}
Lg=&-\frac{1}{M}\big(Q(M,gM)+Q(gM,M)\big),\,
\Gamma(g,g)=\frac{1}{M}Q(gM,gM).
\end{aligned}
\end{equation}

We consider the cases of $0<r<1,0<q$ and $0<q<1,\,0<r$ in this paper. Notice that
\cite{GFSRL05}
\begin{equation}\label{Ker of L}
Ker L=\text{Span}\{1,v_1,v_2,v_3,|v|^2\}.
\end{equation}

\subsection{Well-posedness and hydrodynamic limits}

There has been tremendous progress on the well-posedness of Boltzmann equation. In the context of weak solutions, DiPerna and Lions established the renormalized solutions of the Cauchy problem of Boltzmann equation with large initial datum under Grad's cutoff assumption \cite{DRLPL89}. In \cite{ARVC02}, Alexandre and Villani proved the global existence of renormalized solutions for long-range interaction kernels. For the initial boundary problem, Mischler proved the Boltzmann equation with Maxwell reflection boundary condition for the cutoff case \cite{MS10}. In the context of classical solutions, the global-in-time close to equilibrium classical solution result was first obtained by Ukai in \cite{US74} for collision kernels with cut-off hard potentials. By using the nonlinear energy method, the same type of result for soft potential case for both a periodic domain and for whole space were proved by Guo \cite{GY03,GY04}. The existence and regularity of global classical solution near the equilibrium for whole space were obtained by Gressman-Strain \cite{GPTSRM11} and Alexandre-Morimoto-Ukai-Xu-Yang \cite{ARMYUSXCJYT12} without Grad's angular cutoff kernel assumption. In recent years, there has been significant progress on the strong solutions and mild solutions of the Boltzmann equations in bounded domain endowed with different boundary conditions and were considered in \cite{BMGY16,DRJLWXLLQ22,DRJLSQSSSRM21,GY10,GYKCTDTA17,KCLD18}. Further results references \cite{DRJWY19, ICSL20,ICSLE22,KCLDY18}.

One of the most important features of the Boltzmann equation is that it is connected to fluid equations when the dimensionless number go to zero. Hydrodynamic limits of kinetic equations have been an active research field since the late 1970.

In the context of weak solutions, based on the existence of renormalized solutions \cite{ARVC02, DRLPL89}, Bardos-Golse-Levermore formally derived three fluid equations, including compressible Euler equations, incompressible Euler equations and Navier-Stokes equations, from the scaled Boltzmann equation \cite{BCGFLD91}. They also initialed the program to justify Leray's solutions of the incompressible Navier-Stokes equations from renormalized solutions \cite{CBFGCDL93} under some technical assumptions. This program was completed by Golse and Saint-Raymond with a cutoff Maxwell collision kernel \cite{GFSRL04} and hard cutoff potentials \cite{GFSRL09}. A similar result was obtained by Arsenio for the non-cutoff case \cite{AD12}. Further results references \cite{JNLCDMN10,JNMN17, HLBJJC17, HLBJJC23,MNSRL03}.

In the context of classical solutions, Bardos-Ukai first proved the global existence
of classical solutions $g_\v$ uniformly in $0<\v<1$ and the Navier-Stokes-Fourier limit for cut-off hard potentials \cite{BCUS91}. Employing the semigroup approach, Briant proved incompressible Navier-Stokes limit for Grad's angular cutoff assumption on the torus \cite{BM15}. By applying the Hilbert/Chapman-Enskog expansion, Caflisch, Kawashima-Matsumura-Nishida, Nishida proved the compressible Euler limit, incompressible Navier-Stokes limit from the scaled Boltzmann equation \cite{CRE80,KSMANT79,NT78}. Combining the Hilbert/Chapman-Enskog
expansion with the nonlinear energy method, Guo, Jang and Jiang obtained the acoustic limit \cite{GYJJJN09,GYJJJN10,JJJN09}. Further results references \cite{DMAERLJL89,GY06,GYHFMWY21,JJKC21,JNXLJ15,JNXCJZHJ18}.

For the stationary case, Esposito-Lebowitz-Marra obtained the stationary solution of the stationary Boltzmann equation in a slab with a constant external force and proved that there exists a solution of the scaled Boltzmann equation converges the corresponding compressible Navier-Stokes equations with no-slip boundary conditions \cite{ERLJLMR94}. Esposito-Guo-Kim-Marra applied quantitative $L^2-L^\infty$ approach with new $L^6$ estimates to rigorously derive
the stationary incompressible Navier–Stokes–Fourier system from the stationary Boltzmann equation with a small external field and a small boundary temperature variation for the diffuse boundary condition \cite{ERGYKCMR18}. Wu derived the stationary incompressible Navier–Stokes–Fourier system from the scaled stationary Boltzmann equation in a two dimensional unit plate with in-flow boundary and showed different boundary layer expansion with geometric correction \cite{W16}. Esposito-Marra employed nonlinear estimate to formally obtain the stationary incompressible Navier-Stokes-Fourier limit with Dirichlet boundary conditions from the scaled stationary Boltzmann equation as Knudsen number go to zero on bounded domain or exterior domain \cite{ERMR20}. Further results references \cite{AKGFKS15,GFPF89}.

To our best knowledge, there have no results about derivation of hydrodynamic limits from the scaled Boltzmann equation \eqref{Scaled BE equaion-1} to the stationary equations, including formal derivation and rigorous proofs. In this paper, we formally derive the stationary Navier-Stokes-Fourier system, the stationary Euler equations and the stationary Stokes equations in the Section \ref{Formal limit}.

\section{Notations, difficulty and ideas}
\subsection{Notations}
Before stating main results, we would like to give some notations in this paper.

Denote the spaces $L^2(Mdv)$, as follows:
\begin{equation*}
\begin{aligned}
&L^2(Mdv)=\{f(v)|\int_{\R^3}|f(v)|^2M(v)dv<\infty\},\\
&\l f,g\r_{M}=\int_{\R^3}fgM dv.
\end{aligned}
\end{equation*}

Denote $A=(A_i)$ and $B=(B_{ij})$ the following vectors and tensors:
\begin{equation}\label{The form of A and B}
\begin{aligned}
A_i(v)=&\frac{1}{2}(|v|^2-5)v_i,\,B_{ij}(v)=v_iv_j-\frac{1}{3}|v|^2\delta_{ij}.
\end{aligned}
\end{equation}
Notice that $A_i,B_{ij}\in Ker^{\perp}(L)$ under $\l\c,\c \r_M$ norm.

For $A$ and $B$, we have the following property \cite{CBFGCDL93,DLGF94}:
\begin{proposition}\label{A and B}
There exist $\hat{A}=(\hat{A}_i)$ and $\hat{B}=(\hat{B}_{ij})$ uniquely in $Ker^{\perp}(L)$ such that
\begin{equation}
L(\hat{A}_i)=A_i,\,L(\hat{B}_{ij})=B_{ij}.
\end{equation}
Moreover, there exist two scalar positive functions $\a$ and $\b$ such that
\begin{equation}
\hat{A}(v)=\a(|v|)A(v),\,\hat{B}(v)=\b(|v|)B(v).
\end{equation}
Furthermore, we have
\begin{equation}
\begin{aligned}
\l\hat{A}_i,A_j \r_M&=\frac{5}{2}\kappa\delta_{ij},\\
\l\hat{B}_{ij},B_{kl} \r_M&=\nu(\delta_{ik}\delta_{jl}+\delta_{il}\delta_{jk}-\frac{2}{3}\delta_{ij}\delta_{kl}),\\
\end{aligned}
\end{equation}
where positive constants $\nu$ and $\kappa$ are given by
\begin{equation}
\begin{aligned}
\kappa=&\l\a(|v|),A(v)\otimes A(v)\r_M=\frac{2}{15\sqrt{2\pi}}\int_0^\infty \a (r)r^6 \up{e}^{-\frac{r^2}{2}}dr,\\
\nu=&\l\b(|v|),B(v)\otimes B(v)\r_M=\frac{1}{6\sqrt{2\pi}}\int_0^\infty \b(r)(r^2-5)^2 r^4\up{e}^{-\frac{r^2}{2}}dr.
\end{aligned}
\end{equation}
\end{proposition}

For operators $L$ and $\Gamma$, we also have \cite{BCGFLD91}
\begin{proposition}\label{Poposition of L and Q}
The linearized Boltzmann operator $L$ and the bilinear symmetric operator $\Gamma$ have the following properties:

\up{(i)} The linearized Boltzmann operator $L$ is self-adjoint in $L^2(Mdv)$, i.e. for any $f,g\in L^2(M dv)$, we have
\begin{equation}\label{Self-adjoint}
\l Lf,g\r_M=\l f,Lg\r_M.
\end{equation}

\up{(ii)} For any $g\in Ker L$, we have
\begin{equation}\label{Property of Q}
\frac{1}{2}L (g^2)=\Gamma(g,g).
\end{equation}
\end{proposition}

\subsection{Difficulty and ideas}
In order to derive hydrodynamic limits of the scaled Boltzmann equation \eqref{Scaled BE equaion-1}, the key point is to estimate the singular term $\frac{1}{\v^a}v\c\n_x g_\v$ for some $0<a<1$.

In fact, the scaled Boltzmann equation \eqref{Scaled BE equation} is rewritten as follows:
\begin{equation}\label{Scaled BE equation-3}
\v^{1+q}\p_t g_\v+\v^{q}v\c\n_x g_\v+L g_\v=\v^{r}\Gamma(g_\v,g_\v).
\end{equation}

Suppose $g_\v\to g$ for some function $g$, then formally $\v^{1+q}\p_t g_\v,\,\v^{q}v\c\n_x g_\v,\,\v^{r}\Gamma(g_\v,g_\v)\to\,0$ and $L g_\v\to\, Lg$, therefore $Lg=0$, to combine the fact $KerL=\up{Span}\{M,v_1M,v_2M,v_3M,|v|^2M\}$, we have
\begin{equation}\label{Limits of g}
g(t,x,v)=\rho(t,x)+v\c u(t,x)+(\frac{|v|^2}{2}-\frac{3}{2})\theta(t,x)
\end{equation}
for some functions $\rho,\,u,\,\theta$.

Next, we multiply the scaled Boltzmann equation \eqref{Scaled BE equaion-1} with $vM,\,(\frac{|v|^2}{2}-\frac{5}{2})M$ and integrate by parts over $v\in\R^3$ to get that
\begin{equation}\label{Formally derivation}
\begin{cases}
\v\p_t\l g_\v,v\r_M+\up{div}_{x}\l g_\v,v\otimes v\r_M=0,\\
\v\p_t\l g_\v,\frac{|v|^2}{2}-\frac{5}{2}\r_M+\up{div}_{x}\l g_\v,v(\frac{|v|^2}{2}-\frac{5}{2})\r_M=0.
\end{cases}
\end{equation}

Choosing $\eqref{Formally derivation}_1$ as an example, if $0<r=q<1$, we need to estimate $\frac{1}{\v^r}\up{div}_{x}\l g_\v,v\otimes v\r_M$. Recalling $B_{ij}(v)=v_iv_j-\frac{1}{3}|v|^2\delta_{ij}$ in \eqref{The form of A and B}, we get that
\begin{equation}\label{Key point-1}
\begin{aligned}
\frac{1}{\v^r}\up{div}_{x}\l g_\v,v\otimes v\r_M=&\up{div}_x\l g_\v,B\r_M+\frac{1}{3\v^r}\n_x\l g_\v,|v|^2\r_M\\
=&\frac{1}{\v^r}\up{div}_x\l g_\v,L\hat{B}\r_M+\frac{1}{3\v^r}\n_x\l g_\v,|v|^2\r_M\\
=&\frac{1}{\v^r}\up{div}_x\l Lg_\v,\hat{B}\r_M+\frac{1}{3\v^r}\n_x\l g_\v,|v|^2\r_M,
\end{aligned}
\end{equation}
where Proposition \ref{A and B} is used in the second equality and self-adjoint of $L$ is used in the last equality.

Since $\frac{1}{\v^r}Lg_\v=Q(g_\v,g_\v)-\v\p_t g_\v-v\c\n_x g_\v$ from the scaled Boltzmann equation \eqref{Scaled BE equation}, one of the term in \eqref{Key point-1} is estimated as
\begin{equation}\label{Key point-2}
\begin{aligned}
\frac{1}{\v^r}\up{div}_x\l Lg_\v,\hat{B}\r_M=&\up{div}_x\l Q(g_\v,g_\v),\hat{B}\r_M-\v\up{div}_x\l \p_t g_\v,\hat{B}\r_M-\up{div}_x\l v\c\n_x g_\v,\hat{B}\r_M\\
\,\to&\,\up{div}_x\l Q(g,g),\hat{B}\r_M-\up{div}_x\l v\c\n_x g,\hat{B}\r_M\\
=&\frac{1}{2}\up{div}_x\l L(g^2),\hat{B}\r_M-\nu\Delta_x u\\
=&\frac{1}{2}\up{div}_x\l g^2,B\r_M-\nu\Delta_x u\\
=&u\c\n_x u-\nu\Delta_x u,
\end{aligned}
\end{equation}
where we use Proposition \ref{A and B}, self-adjoint of $L$ and the form of g in \eqref{Limits of g}.

To deal with the singular term $\frac{1}{3\v^r}\n_x\l g_\v,|v|^2\r_M$, applying the Leray projection $\mathbb{P}$ on $\eqref{Key point-1}_1$, then we only estimate
\begin{equation}\label{Key point-3}
\begin{aligned}
\frac{1}{\v^r}\mathbb{P}\up{div}_x\l g_\v,v\otimes v\r_M=\frac{1}{\v^r}\up{div}_x\l g_\v,B\r_M.
\end{aligned}
\end{equation}

Therefore, we formally derive the stationary Navier-Stokes-Fourier system (see Theorem \ref{Main theorem}):
\begin{equation*}
\begin{cases}
u\c\n_x u+\n_x p=\nu \Delta_x u,\\
\up{div}_x u=0,\\
\n_x(\rho+\theta)=0,\\
u\c\n_x\theta=\kappa\Delta_x\theta.
\end{cases}
\end{equation*}

\section{Main results and the proof} \label{Formal limit}

\subsection{Main result}
The main theorem is stated as follows:
\begin{theorem}\label{Main theorem}
Let $f_\v(t, x, v)$ be a sequence of nonnegative solutions to the scaled Boltzmann equation \eqref{Scaled BE equaion-1} with a formula of \eqref{Special form solution}, the sequence $g_\v$ converges to a function $g_\v$ as $\v$ goes to zero in the sense of distributions. Furthermore, assume the moments that
\begin{equation}
\begin{aligned}
&\l g_\v, 1\r_M,\,\l g_\v, v\r_M,\,\l g_\v,v\otimes v\r_M,\,\l g_\v,v|v|^2\r_M,\\
&\l g_\v, \hat{A}(v)\otimes v\r_M,\,\l \Gamma(g_\v,g_\v),\hat{A}(v)\r_M,\\
&\l g_\v, \hat{B}(v)\otimes v\r_M,\,\l \Gamma(g_\v,g_\v),\hat{B}(v)\r_M
\end{aligned}
\end{equation}
converge in $D'(\R\times\Omega)$ to the corresponding moments
\begin{equation}
\begin{aligned}
&\l g,1\r_M,\,\l v,g\r_M,\,\l g,v\otimes v\r_M,\,\l g,v|v|^2\r_M,\\
&\l g,\hat{A}(v)\otimes v\r_M,\,\l \Gamma(g,g), \hat{A}(v)\r_M,\\
&\l g,\hat{B}(v)\otimes v\r_M,\,\l \Gamma(g,g),\hat{B}(v)\r_M
\end{aligned}
\end{equation}
as $\v\to 0^+$. Then the limiting $g$ has the form
\begin{equation}\label{The form of g}
g=\rho+u\c v+\frac{1}{2}(|v|^2-3)\theta,
\end{equation}
where the velocity $u$ is divergence-free and the density and temperature fluctuations $\rho$ and $\theta$ satisfy the Boussinesq relation
\begin{equation}\label{The Boussinesq relation}
\up{div}_xu=0,\,\n_x(\rho+\theta)=0.
\end{equation}
Moreover, the functions $\rho,\,u$ and $\theta$ are solutions of the equations
\begin{equation}\label{The limits of Eqs}
\begin{aligned}
u\c\n_x u+\n_x p=&\nu\Delta_x u,\,\,\,\,\quad u\c\n_x\theta=\kappa\Delta_x\theta,\qquad\qquad\qquad\up{if}\,0<r=q<1,\\
\n_x p=&\nu\Delta_x u,\quad\quad\,\,\kappa\Delta_x\theta=0,\qquad\qquad\qquad\up{if}\,0<q<\min\{1,r\},0<r,\\
u\c\n_x u+\n_x p=&0,\,\qquad\quad\,\,u\c\n_x\theta=0,\qquad\qquad\qquad\up{if}\,0<r<\min\{1,q\},0<q,
\end{aligned}
\end{equation}
where $\nu$ and $\kappa$ are given in Proposition \ref{A and B}.
\end{theorem}
\begin{remark}
In \cite{BCGFLD91}, Bardos, Glose and Levermore formally derived the Navier-Stokes-Fourier system and other three equations by considering $r,q\geq1$. It's different from non-stationary limit if we want to give a rigorous proof from the scaled Boltzmann equation to the stationary equations.

Taking the Navier-Stokes-Fourier limit as an example. On the one hand, for the non-stationary Navier-Stokes-Fourier limit, we have the following scaled Boltzmann equation
\begin{equation}\label{The Navier-Stokes-Fourier limit}
\p_t g_\v+\frac{1}{\v}v\c\n_x g_\v+\frac{1}{\v^2}Lg_\v=\frac{1}{\v}\Gamma (g_\v,g_\v).
\end{equation}
While following scaled Boltzmann equation shows the stationary Navier-Stokes-Fourier limit
\begin{equation}\label{The stationary Navier-Stokes-Fourier limit}
\p_t g_\v+\frac{1}{\v}v\c\n_x g_\v+\frac{1}{\v^{1+r}}Lg_\v=\frac{1}{\v}\Gamma(g_\v,g_\v).
\end{equation}

Formally, the dissipation term is $O(\frac{1}{\v^2})$ from the linear Boltzmann operator $L$ in \eqref{The Navier-Stokes-Fourier limit}. However, the order of the dissipation term is $\frac{1}{\v^{1+r}}$ in The stationary Navier-Stokes-Fourier limit, which is weaker than $\frac{1}{\v^2}$ since $0<r<1$. Furthermore, at least in the context of classical solutions, to prove $\p_t\l g_\v,v\r_M$ is $O(\v^a)$ uniformly with $\v$ for some positive constant $a>0$ if we want to give a rigorous proof for the stationary case, while we just need to prove the term $\p_t\l g_\v,v\r_M$ is $O(1)$ for the non-stationary case, which means the proof of the stationary case is maybe more difficult in this sense. Therefore, we need to looking for additional dissipative structures to prove the stationary Navier-Stokes-Fourier limit.

On the other hand, for the scaled Boltzmann equation, initial conditions need to be provided, while the stationary equations for the limiting equations does not require initial conditions. This implies that there will be an initial boundary layer. In this sense, the stationary equations limit is more complex than the non-stationary equations.
\end{remark}

\begin{proof}

Recalling the scaled Boltzmann equation \eqref{Scaled BE equation} as follow:
\begin{equation}\label{Scaled BE equation-1}
\v\p_t g_\v+v\c\n_x g_\v+\frac{1}{\v^{q}}L g_\v=\frac{1}{\v^{q-r}}\Gamma(g_\v,g_\v).
\end{equation}
The scaled Boltzmann equation \eqref{Scaled BE equation-1} is rewritten as:
\begin{equation}\label{Scaled BE equation-2}
L g_\v=-\v^{1+q}\p_t g_\v-\v^ qv\c\n_x g_\v-\v^{r}\Gamma(g_\v,g_\v).
\end{equation}

Letting $\v\to 0^+$ and using the assumption of moment convergence implies the relation
\begin{equation}\label{Lg}
L g=0
\end{equation}
which shows that $g$ belongs to $Ker L$ \eqref{Ker of L} and can be written in the form of \eqref{The form of g}. 

The derivation of \eqref{The Boussinesq relation} from the following conservation of mass and momentum:
\begin{equation}
\begin{aligned}
\v\p_t\l g_\v, 1\r_M+\up{div}_x\l g_\v, v\r_M&=0,\\
\v\p_t\l g_\v, v\r_M+\up{div}_x\l g_\v ,v\otimes v\r_M&=0.
\end{aligned}
\end{equation}
Letting $\v\to 0^+$ in the above expression, we can obtain the relations
\begin{equation}
\begin{aligned}
\up{div}_x\l g,v\r_M&=0,\,\up{div}_x\l g ,v\otimes v\r_M=0.
\end{aligned}
\end{equation}
Substituting $g$ in \eqref{The form of g} into the left-hand side of the above relations, we can obtain  \eqref{The Boussinesq relation}. Performing the same operation on the conservation of energy also yields the divergence-free velocity condition of \eqref{The Boussinesq relation}.

If $0<r=q<1$, the scaled Boltzmann equation \eqref{Scaled BE equation-1} is rewritten as
\begin{equation}\label{Scaled BE equation-3}
\v^{1-r}\p_t g_\v+\frac{1}{\v^r}v\c\n_x g_\v+\frac{1}{\v^{2r}}L g_\v=\frac{1}{\v^r}\Gamma(g_\v,g_\v).
\end{equation}

Then, we deduce that
\begin{equation}\label{The case of r=q-1}
\begin{aligned}
\v^{1-r}\p_t\l g_\v, v\r_M+\frac{1}{\v^r}\up{div}_x\l g_\v,v\otimes v\r_M&=0,\\
\v^{1-r}\p_t\l g_\v,(\frac{|v|^2}{2}-\frac{5}{2})\r_M+\frac{1}{\v^r}\up{div}_x\l g_\v ,(\frac{|v|^2}{2}-\frac{5}{2})v\r_M&=0.
\end{aligned}
\end{equation}

Using of the moment convergence assumption, the form of $g$ given by \eqref{The form of g} and $0<r<1$ provides
\begin{equation}\label{The case of limits r=q-1}
\begin{aligned}
&\lim_{\v\to 0^+}\v^{1-r}\p_t\l g_\v, v\r_M=0,\\
&\lim_{\v\to 0^+}\v^{1-r}\p_t\l g_\v,(\frac{|v|^2}{2}-\frac{5}{2})\r_M=0.
\end{aligned}
\end{equation}

To complete the proof of $\eqref{The limits of Eqs}_1$, it is need to estimate $\v^{-r}\up{div}_x\l g_\v,v\otimes v\r_M$ and $\v^{-r}\up{div}_x\l g_\v,(\frac{|v|^2}{2}-\frac{5}{2})v\r_M$ in \eqref{The case of r=q-1}.

Applying the self-adjoint of the linear operator $L$  in Poposition \ref{Poposition of L and Q} and Proposition \ref{A and B}, we have
\begin{equation}\label{The case of limits r=q-2}
\begin{aligned}
\lim_{\v\to 0^+}\mathbb{P}\frac{1}{\v^r}\up{div}_x\l g_\v,v\otimes v\r_M=&\lim_{\v\to 0^+}\frac{1}{\v^r}\up{div}_x\l g_\v,B\r_M=\lim_{\v\to 0^+}\frac{1}{\v^r}\up{div}_x\l g_\v,L\hat{B}\r_M\\
=&\lim_{\v\to 0^+}\frac{1}{\v^r}\up{div}_x\l Lg_\v,\hat{B}\r_M,\\
\lim_{\v\to 0^+}\frac{1}{\v^r}\up{div}_x\l g_\v,(\frac{|v|^2}{2}-\frac{5}{2})v\r_M=&\lim_{\v\to 0^+}\frac{1}{\v^r}\up{div}_x\l g_\v,A\r_M=\lim_{\v\to 0^+}\frac{1}{\v^r}\up{div}_x\l g_\v,L\hat{A}\r_M\\
=&\lim_{\v\to 0^+}\frac{1}{\v^r}\up{div}_x\l Lg_\v,\hat{A}\r_M,
\end{aligned}
\end{equation}
where the symbol $\mathbb{P}$ is the Leray projection.

Recalling the scaled Boltzmann equation in \eqref{Scaled BE equation-3}
\begin{equation}\label{Scaled BE equation-4}
\frac{1}{\v^r}L g_\v=\Gamma(g_\v,g_\v)-\v\p_t g_\v-v\c\n_x g_\v.
\end{equation}

Then, to apply the fact that $\frac{1}{2}L(g^2)=\Gamma(g,g)$, we have
\begin{equation*}
\begin{aligned}
\lim_{\v\to 0^+}\mathbb{P}\frac{1}{\v^r}\up{div}_x\l g_\v,v\otimes v\r_M=&\lim_{\v\to 0^+}\up{div}_x\l \Gamma(g_\v,g_\v),\hat{B}\r_M-\lim_{\v\to 0^+}\v\up{div}_x\l \p_t g_\v,\hat{B}\r_M\\
&-\lim_{\v\to 0^+}\up{div}_x\l v\c\n_x g_\v,\hat{B}\r_M\\
=&\frac{1}{2}\up{div}_x\l L(g^2),\hat{B}\r_M-\up{div}_x\l v\c\n_x g,\hat{B}\r_M\\
=&\frac{1}{2}\up{div}_x\l g^2,L\hat{B}\r_M-\up{div}_x\l v\c\n_x g,\beta(|v|)B\r_M\\
=&u\c\n_x u-\nu\Delta_x u,\\
\end{aligned}
\end{equation*}
and
\begin{equation*}
\begin{aligned}
\lim_{\v\to 0^+}\frac{1}{\v^r}\up{div}_x\l g_\v,(\frac{|v|^2}{2}-\frac{5}{2})v\r_M=&\lim_{\v\to 0^+}\up{div}_x\l \Gamma(g_\v,g_\v),\hat{A}\r_M-\lim_{\v\to 0^+}\v\up{div}_x\l \p_t g_\v,\hat{A}\r_M\\
&-\lim_{\v\to 0^+}\up{div}_x\l v\c\n_x g_\v,\hat{A}\r_M\\
=&\frac{1}{2}\up{div}_x\l L(g^2),\hat{A}\r_M-\up{div}_x\l v\c\n_x g,\hat{A}\r_M\\
=&\frac{1}{2}\up{div}_x\l g^2,L\hat{A}\r_M-\up{div}_x\l v\c\n_x g,\alpha(|v|)A\r_M\\
=&u\c\n_x\theta-\kappa\Delta_x \theta,
\end{aligned}
\end{equation*}
where Proposition \ref{A and B} and $g=\rho+u\c v+(\frac{|v|^2}{2}-\frac{3}{2})\theta$ in \eqref{The form of g} are used.

Therefore, the stationary incompressible Navier-Stokes-Fourier system is derived as follows:
\begin{equation}\label{Stationary Navier-Stokes Eq}
\begin{cases}
u\c\n_xu+\n_x p=\nu\Delta_x u,\\
\up{div}_x u=0,\\
u\c\n_x\theta=\kappa\Delta_x\theta,\\
\n_x(\rho+\theta)=0.
\end{cases}
\end{equation}

If $0<q<\min\{1,r\},\,0<r$ or $0<r<\min\{1,q\},\,0<q$, the scaled Boltzmann equation \eqref{Scaled BE equation-1} is rewritten as
\begin{equation}\label{Scaled BE equation-4}
\v^{1-q}\p_t g^\v+\frac{1}{\v^q}v\c\n_x g_\v+\frac{1}{\v^{2q}}L g_\v=\v^{r-2q}\Gamma(g_\v,g_\v),
\end{equation}
and
\begin{equation}\label{Scaled BE equation-5}
\v^{1-r}\p_t g_\v+\frac{1}{\v^r}v\c\n_x g_\v+\frac{1}{\v^{r+q}}L g_\v=\frac{1}{\v^q}\Gamma(g_\v,g_\v).
\end{equation}

Then, applying a similar estimate as the case of $0<r=q<1$, we can derive the equations of $\eqref{The limits of Eqs}_2$ and $\eqref{The limits of Eqs}_3$.

\end{proof}


%


%
%
%
%
%
%
%
%
%
%


\end{document}